\documentclass{amsart}

\usepackage[scale=0.76]{geometry}

\usepackage{mathtools}
\usepackage{amssymb, hyperref}
\usepackage{amsmath,amssymb,amsbsy,amsfonts, latexsym,amsopn,amstext, amsxtra,euscript,amscd}
\usepackage{amsthm}

\usepackage{rotating}

\usepackage{amssymb, amsthm, amsmath,  amssymb, amsbsy,   amsfonts,  latexsym,  amsopn,   amstext, amsxtra,  euscript,   amscd}
\usepackage{lscape} 

\usepackage{cite}
\usepackage{color}
\usepackage[shortlabels]{enumitem}

\usepackage{tikz}
\usetikzlibrary{decorations.pathreplacing,decorations.markings}

\usetikzlibrary{positioning}

\usepackage[capitalise]{cleveref}

 \usepackage[all,poly]{xy}

\usepackage[alphabetic, citation-order, msc-links]{amsrefs}

\hypersetup{
  colorlinks   = true, 
  urlcolor     = blue, 
  linkcolor    = blue, 
  citecolor   = red 
}

\newtheorem{thm}{Teorema}
\newtheorem{prop}{Pohim}
\newtheorem{lem}{Lema}

\newtheorem{rem}{Shënim}
\newtheorem{exa}{Shembull}

\theoremstyle{definition}

\theoremstyle{remark}

\crefname{thm}{Thm.}{}
\crefname{prop}{Prop.}{}
\crefname{lem}{Lem.}{}
\crefname{cor}{Cor.}{}
\DeclareMathOperator\gl{GL}     	 
\DeclareMathOperator\Sl{SL} 
\newcommand{\card}[1]{\left\lvert\mspace{1mu}#1\mspace{1mu}\right\rvert}

\def\X{\mathcal X}

\def\O{\mathcal O}
\def\p{\mathfrak p}

\def\s{\mathfrak s}

\DeclareMathOperator\Stab{Stab }
\DeclareMathOperator\chara{char }

\def\F{\mathbb F}

\def\X{\mathcal X}

\def\iso{\cong}

\def\g{\gamma}

\def\t{\tau}
\def\d{\delta}
\def\a{\alpha}
\def\b{\beta}

\def\D{\Delta}

\DeclareMathOperator\Orb{Orb}

\def\l{\lambda}

\DeclareMathOperator\Gal{Gal}

\def\s{\sigma}
\def\<{\left<}
\def\>{\right>}
\def\es{është }

\begin{document}

\title{Zgjidhja e ekuacionit të gradës 5-të}
\author{Elira Shaska}
\thanks{Supported by $A^2G^2$-group at RISAT}


\date{March 15, 2021}

\begin{abstract} 
Ekuacioni polinomial $f(x)=0$   ka zgjidhje me radikale atëherë dhe vetëm atëherë kur grupi i Galuait $\Gal_\F (f)$  është grup i zgjidhshëm.  Ne japim kushte të nevojshme dhe të mjaftueshme për përcaktimin e $\Gal_\F (f)$ dhe në rastet kur $\Gal_\F (f)$ është grup i zgjidhshëm japim formulat përkatëse të gjetjes së rrënjëve. 
\end{abstract}

\maketitle

\section{Hyrje}

Gjetja e formulave për  rrënjët e  ekuacioneve polinomialë është një nga problemet më klasike të matematikës. Formula kuadratike (për gjetjen e rrënjëve të polinomeve kuadratike) është nga formulat e para të algjebrës.  Për polinome me gradë $n=3$ (kubiket)  këto formula janë disi më të komplikuara dhe njihen si formulat e Kardanos.  Pës $n=4$ (kuartiket)   formulat e gjetjes së rrënjëve quhen formulat e Ferrarit; shih \cite{algjebra} për më tepër detaje.  

Ekuacionet e gradës $n\geq 5$ (kuintiket) paraqitën vështiresi më të mëdha.  Në 1799 Ruffini paraqiti një "vërtetim të pasaktë" të faktit se nuk ekziston një formulë për zgjidhjen e ekuacionit kuintik.   Në 1824, Nils  Abel ishte i pari  që dha një vërtetim të saktë  të faktit që      një formulë e gjetjes së rrënjëve për kuintiket nuk ekziston. Këtij fakti tashmë në literaturë i referohet si teorema Abel-Ruffini; shihni  \cite{king-beyond-quartic}. 

Me zhvillimin e Teorisë së Galuait, teoria a gjetjes së formulave që përcaktojnë rrënjët (zgjidhjet me radikale) u transformua dhe u përgjithësua për polinomet me gradë çfarëdo. Rezultati kryesor i Teorisë së Galuait \es: \emph{një ekuacion polinomial $f(x)=0$, për $f\in \F[x]$, është  i zgjidhshëm me radikale atëherë dhe vetëm atëherë kur $\Gal_\F (f)$ është një grup i zgjidhshëm.}  Natyrshëm lind pyetja: \emph{cilat janë formulat për ekuacionet kuintike që janë të zgjidhshme?}

Në këtë përmbledhje të disa fakteve bazë mbi teorinë e zgjidhjeve të ekuacioneve ne synojmë të japim një përshkrim të shkurtër të zgjidhjes së ekuacioneve kuintike.  
Ne supozojmë se lexuesi është i familiarizuar me faktet bazë të algjebrës në nivelin e \cite{algjebra}

\section{Njohuri paraprake}

Le të jetë $\F$ një fushë perfekte me karakteristike $\chara \F=0$.   Jepet  $f(x)\in \F[x]$ një polinom me gradë $n$ dhe  $\a_1, \dots, \a_n$  rrënjët e këtij polinomi. Fusha   $\F(\a_1, \dots, \a_n)$ quhet  \textbf{fusha ndarëse}  e polinomit $f(x)$. Nga një teoremë e Kronecker, për çdo fushë $\F$ dhe çdo polinom  $f(x)\in \F[x]$  ekziston një shtrirje  $E_f$ e $\F$-së që është fusha ndarëse e   $f(x)$.  

\begin{thm}[Teorema e shtrirjes së izomorfizmit] \label{iso-ext-thm}
Le të jetë $\s: \F \to \F^\prime$ një izomorfizëm i fushës. Kemi $S=\{ f_i (x) \}$ një bashkësi polinomesh në $\F$ dhe $S^\prime=\{ \s (f_i ) \}$ polinomet korresponduese në $\F^\prime$. Le të jetë $K$ fusha ndarëse për $S$ mbi $\F$  dhe  $K^\prime$ fusha ndarëse $S^\prime$ mbi $\F'$. Atëherë ekziston një automorfizëm $\t : K \to K^\prime$ me $\t \mid_\F=\s$.

Për më shumë, nëse $\a \in K$ dhe $\a^\prime$ është një rrënjë e polinomit $\s (\min (\a, \F, x))$, atëherë  $\t$ mund të zgjidhet i tillë që $\t (\a) = \a^\prime$.
\end{thm}

 Pra, për çdo  $f(x) \in \F[x]$ çdo dy fusha ndarëse të  $f(x)$-së janë   $\F$-izomorfike.

\begin{lem}\label{sp_field} 
Fusha ndarëse e një polinomi me gradë $n$ mbi $\F$ ka gradë $[E_f : \F]\leq n!$. Në qoftë se $f(x)$ është i pathjeshtueshëm mbi $\F$, atëherë  $n \mid   [E_f:\F]$.
\end{lem}

\begin{thm}\label{thm-13-2-1}
Le të jetë $ f \in \F[x]$ dhe    $E_f$ fusha ndarëse mbi $\F$. Nëse $ f(x)$ nuk ka rrënjë të shumëfishta, atëherë 
\begin{equation} 
\card{\Gal (E_f/ \F) } = [E_f :F]. 
\end{equation} 
\end{thm}

\subsection{Grupi Galua i polinomeve}
Jepet polinomi $ f(x) $ me gradë $n=\deg f$ , i pafaktorizueshëm në    $ \F[x]$, i cili   faktorizohet si
\begin{equation} f(x)= (x-\a_1) \dots (x-\a_n)\end{equation} 
në fushën ndarëse $E_f$. Ateherë,  $E_f/\F$ është shtrirje  Galua sepse është shtrirje normale ($E_f$ është fushë ndarëse) dhe shtrirje e ndashme (sepse $\F$ është fushë perfekte). Grupi    $ \Gal (E_f/\F)$ quhet  \textbf{grupi i Galuait} i $ f(x)$ mbi $\F$ dhe shënohet me  $\Gal_\F (f)$. Elementët e  $\Gal_\F (f)$  përkëmbejnë rrënjët e  $ f(x)$. Pra grupi i Galuait i një polinomi ka një kopje izomorfike në $S_n$, e cila përcaktohet deri në klasën e konjugimit nga $f(x)$

\begin{prop} Dy pohimet e mëposhtme janë të vërteta:
\begin{enumerate}[\upshape(i), nolistsep] 
\item Shënojmë me  $G=\Gal_\F (f)$ dhe $H=G\cap A_n$. Atëherë $H= \Gal (E_f / \F(\sqrt{\D_f}))$. Për më tepër,   $G$ është një nëngrup i grupit alternativ $A_n$ atëherë dhe vetëm atëherë kur   diskriminanti  $\D_f$ është një katror i plotë në  $\F$.
\item Faktorët e pazbërthyeshëm të  $f$ në  $\F[x]$ korrespondojnë me orbitat e $G$-së. Për më tepër,  $G$ është një nëngrup transitiv i  $S_n$ atëherë dhe vetëm atëherë kur   $f$ është i pathjeshtueshëm. 
\end{enumerate} 
\end{prop}

Për vërtetimin e këtyre fakteve shihni \cite{algjebra}. 

\section{Zgjidhja e ekuacioneve polinomialë të gradës së 5-të}
\begin{lem}
Le të jetë $f(x)\in \F[x]$ një polinom i pathjeshtueshëm i gradës së pestë. Atëherë $\Gal_\F (f) $ është izomorfik  me një nga këta grupe: $C_5, D_5$,  $F_5= AGL(1,5)$, $A_5$, or $S_5$
\end{lem}

\proof $G$ është transitiv, prandaj nëngrupi i tij Sylow-5 është izomorfik me  $C_5$ (i gjeneruar nga një cikël i rendit 5). Nëse $C_5$ është normal, atëherë $G$ ka të paktën 6 nëngrupe Sylow-5; atëherë $\card{G}\ge 6\cdot 5=30$, prandaj $[S_5:G]\le 4$ nga ku rrjedh $G = S_5, A_5$. Nëse $C_5$ është normal në $G$ atëherë $G$ është i  konjuguar ose  $C_5$, $D_5$ (grupi dihedral i rendit 10) ose $F_5= AGL(1,5)$, normalizuesi i plotë i  $C_5$ në $S_5$, i rendit 20 (i quajtur edhe grupi Frobenius i rendit 20).
\qed

\begin{rem}
Nëse dallori i polinomit të gradës së pestë është katror i plotë në $\F$ atëherë $\Gal (f) $ përmbahet në $A_5$. Kjo do të thotë që është  $C_5, D_5$, ose $A_5$.
\end{rem}

\subsection{Kuintikët e zgjidhshëm}

Nëse $ G = S_5,  A_5$ atëherë ekuacioni $f(x) = 0$ nuk është i zgjidhshëm me radikale. Ne duam të shqyrtojmë këtu rastin $G \not\iso  S_5, \ A_5$.

Le të jetë $f(x)$ një polinom i pathjeshtueshëm i gradës së pestë në $\F[x]$ i dhënë në formën:
\[
f(x) =  x^5\ +\ c_4x^4\ +   \cdots \ + \ c_0   =  (x-\a_1) \cdots  (x-\a_5)\leqno{(5)}
\]
Le të jetë $G = \Gal (f)$, nëngrup transitiv  $S_5$ i përftuar nga përkëmbimi i rrënjëve të ndryshme $\a_1, \cdots,\a_5$. Si më parë  $E_f = k(\a_1, \cdots,\a_5)$ përfaqëson fushën ndarëse të $f$-së. Një element 5-ciklor në  $S_5=\mbox{Sym}(\{1, \dots,5\})$ i korrespondon një pentagoni (të orientuar)  me kulme  $1, \dots ,5$.  Një  5-cikël dhe inversi i tij korrespondojnë me një pentagon të pa orientuar.

\medskip

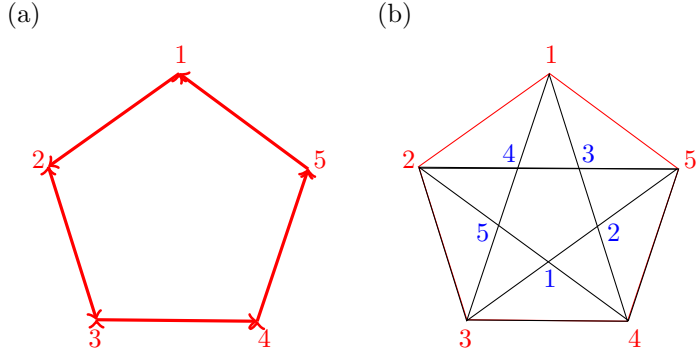
\begin{tabular}{cc}
\begin{minipage}{1.8in}
(a)\begin{displaymath}
\begin{tikzpicture}[x=0.50pt,y=0.50pt,yscale=-1, xscale=1] 
\draw[->, red, very thick]  (95.38,32.87) edge (-3.6,103.37);
\draw[->, red, very thick]  (-3.6,103.37) edge  (32.86,219.29);
\draw[->, red, very thick]  (32.86,219.29) edge (154.38,220.43);
\draw[->, red, very thick] (154.38,220.43) edge (193.02,105.22) ;
\draw[->, red, very thick] (193.02,105.22)  edge (95.38,32.87);
\draw[red] (90,10) node [anchor=north west][inner sep=0.75pt]   [align=left] {1};
\draw[red] (-18,90) node [anchor=north west][inner sep=0.75pt]   [align=left]  {2};
\draw[red] (25,225) node [anchor=north west][inner sep=0.75pt]   [align=left] {3};
\draw[red] (153,225) node [anchor=north west][inner sep=0.75pt]   [align=left] {4};
\draw[red] (195,90) node [anchor=north west][inner sep=0.75pt]   [align=left] {5};
\end{tikzpicture}
\end{displaymath}
\end{minipage}
&
\begin{minipage}{1.8in}
(b)\begin{displaymath}
\begin{tikzpicture}[x=0.50pt,y=0.50pt,yscale=-1,xscale=1]  
\draw[red]   (193.02,105.22) -- (154.38,220.43) -- (32.86,219.29) -- (-3.6,103.37) -- (95.38,32.87) -- cycle ;
\draw   (95.12,33.52) -- (154.38,220.43) -- (32.94,219.48) -- cycle ;
\draw   (-3.26,103.84) -- (193.02,105.22) -- (155.06,219.87) -- cycle ;
\draw   (191.79,105.18) -- (32.71,220.07) -- (-3.3,104.53) -- cycle ;
\draw[red] (90,10) node [anchor=north west][inner sep=0.75pt]   [align=left]  {1};
\draw[red] (-18,90) node [anchor=north west][inner sep=0.75pt]   [align=left]  {2};
\draw[red] (25,225) node [anchor=north west][inner sep=0.75pt]   [align=left]  {3};
\draw[red] (153,225) node [anchor=north west][inner sep=0.75pt]   [align=left]  {4};
\draw[red] (195,90) node [anchor=north west][inner sep=0.75pt]   [align=left]  {5};
\draw[blue] (89,180) node [anchor=north west][inner sep=0.75pt]   [align=left]  {1};
\draw[blue] (137,145) node [anchor=north west][inner sep=0.75pt]   [align=left]  {2};
\draw[blue] (118,85) node [anchor=north west][inner sep=0.75pt]   [align=left]  {3};
\draw[blue] (58 ,85) node [anchor=north west][inner sep=0.75pt]   [align=left]  {4};
\draw[blue] (38,145) node [anchor=north west][inner sep=0.75pt]   [align=left]  {5};
\end{tikzpicture}
\end{displaymath}
\end{minipage}
 \end{tabular}

\medskip

Kurse gjithë grupi ciklik   $C_5$ i korrespondon pentagonit së bashku me "të kundërtin e tij".

\medskip

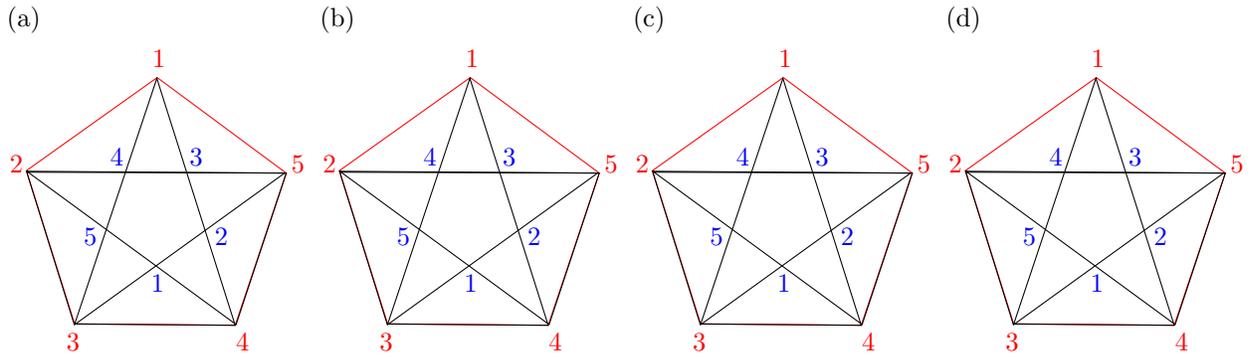
\begin{tabular}{cccc}
\begin{minipage}{1.5in}
(a)\begin{displaymath}
\begin{tikzpicture}[x=0.50pt,y=0.50pt,yscale=-1,xscale=1]  
\draw[red]   (193.02,105.22) -- (154.38,220.43) -- (32.86,219.29) -- (-3.6,103.37) -- (95.38,32.87) -- cycle ;
\draw   (95.12,33.52) -- (154.38,220.43) -- (32.94,219.48) -- cycle ;
\draw   (-3.26,103.84) -- (193.02,105.22) -- (155.06,219.87) -- cycle ;
\draw   (191.79,105.18) -- (32.71,220.07) -- (-3.3,104.53) -- cycle ;
\draw[red] (90,10) node [anchor=north west][inner sep=0.75pt]   [align=left]  {1};
\draw[red] (-18,90) node [anchor=north west][inner sep=0.75pt]   [align=left]  {2};
\draw[red] (25,225) node [anchor=north west][inner sep=0.75pt]   [align=left]  {3};
\draw[red] (153,225) node [anchor=north west][inner sep=0.75pt]   [align=left]  {4};
\draw[red] (195,90) node [anchor=north west][inner sep=0.75pt]   [align=left]  {5};
\draw[blue] (89,180) node [anchor=north west][inner sep=0.75pt]   [align=left]  {1};
\draw[blue] (137,145) node [anchor=north west][inner sep=0.75pt]   [align=left]  {2};
\draw[blue] (118,85) node [anchor=north west][inner sep=0.75pt]   [align=left]  {3};
\draw[blue] (58 ,85) node [anchor=north west][inner sep=0.75pt]   [align=left]  {4};
\draw[blue] (38,145) node [anchor=north west][inner sep=0.75pt]   [align=left]  {5};
\end{tikzpicture}
\end{displaymath}
\end{minipage}
&
\begin{minipage}{1.5in}
(b)\begin{displaymath}
\begin{tikzpicture}[x=0.50pt,y=0.50pt,yscale=-1,xscale=1]  
\draw[red]   (193.02,105.22) -- (154.38,220.43) -- (32.86,219.29) -- (-3.6,103.37) -- (95.38,32.87) -- cycle ;
\draw   (95.12,33.52) -- (154.38,220.43) -- (32.94,219.48) -- cycle ;
\draw   (-3.26,103.84) -- (193.02,105.22) -- (155.06,219.87) -- cycle ;
\draw   (191.79,105.18) -- (32.71,220.07) -- (-3.3,104.53) -- cycle ;
\draw[red] (90,10) node [anchor=north west][inner sep=0.75pt]   [align=left]  {1};
\draw[red] (-18,90) node [anchor=north west][inner sep=0.75pt]   [align=left]  {2};
\draw[red] (25,225) node [anchor=north west][inner sep=0.75pt]   [align=left]  {3};
\draw[red] (153,225) node [anchor=north west][inner sep=0.75pt]   [align=left]  {4};
\draw[red] (195,90) node [anchor=north west][inner sep=0.75pt]   [align=left]  {5};
\draw[blue] (89,180) node [anchor=north west][inner sep=0.75pt]   [align=left]  {1};
\draw[blue] (137,145) node [anchor=north west][inner sep=0.75pt]   [align=left]  {2};
\draw[blue] (118,85) node [anchor=north west][inner sep=0.75pt]   [align=left]  {3};
\draw[blue] (58,85) node [anchor=north west][inner sep=0.75pt]   [align=left]  {4};
\draw[blue] (38,145) node [anchor=north west][inner sep=0.75pt]   [align=left]  {5};
\end{tikzpicture}
\end{displaymath}
\end{minipage}
&
\begin{minipage}{1.5in}
(c)\begin{displaymath}
\begin{tikzpicture}[x=0.50pt,y=0.50pt,yscale=-1,xscale=1]  
\draw[red]   (193.02,105.22) -- (154.38,220.43) -- (32.86,219.29) -- (-3.6,103.37) -- (95.38,32.87) -- cycle ;
\draw   (95.12,33.52) -- (154.38,220.43) -- (32.94,219.48) -- cycle ;
\draw   (-3.26,103.84) -- (193.02,105.22) -- (155.06,219.87) -- cycle ;
\draw   (191.79,105.18) -- (32.71,220.07) -- (-3.3,104.53) -- cycle ;
\draw[red] (90,10) node [anchor=north west][inner sep=0.75pt]   [align=left]  {1};
\draw[red] (-18,90) node [anchor=north west][inner sep=0.75pt]   [align=left]  {2};
\draw[red] (25,225) node [anchor=north west][inner sep=0.75pt]   [align=left]  {3};
\draw[red] (153,225) node [anchor=north west][inner sep=0.75pt]   [align=left]  {4};
\draw[red] (195,90) node [anchor=north west][inner sep=0.75pt]   [align=left]  {5};
\draw[blue] (89,180) node [anchor=north west][inner sep=0.75pt]   [align=left]  {1};
\draw[blue] (137,145) node [anchor=north west][inner sep=0.75pt]   [align=left]  {2};
\draw[blue] (118,85) node [anchor=north west][inner sep=0.75pt]   [align=left]  {3};
\draw[blue] (58,85) node [anchor=north west][inner sep=0.75pt]   [align=left]  {4};
\draw[blue] (38,145) node [anchor=north west][inner sep=0.75pt]   [align=left]  {5};
\end{tikzpicture}
\end{displaymath}
\end{minipage}
&
\begin{minipage}{1.5in}
(d)\begin{displaymath}
\begin{tikzpicture}[x=0.50pt,y=0.50pt,yscale=-1,xscale=1]  
\draw[red]   (193.02,105.22) -- (154.38,220.43) -- (32.86,219.29) -- (-3.6,103.37) -- (95.38,32.87) -- cycle ;
\draw   (95.12,33.52) -- (154.38,220.43) -- (32.94,219.48) -- cycle ;
\draw   (-3.26,103.84) -- (193.02,105.22) -- (155.06,219.87) -- cycle ;
\draw   (191.79,105.18) -- (32.71,220.07) -- (-3.3,104.53) -- cycle ;
\draw[red] (90,10) node [anchor=north west][inner sep=0.75pt]   [align=left]  {1};
\draw[red] (-18,90) node [anchor=north west][inner sep=0.75pt]   [align=left]  {2};
\draw[red] (25,225) node [anchor=north west][inner sep=0.75pt]   [align=left]  {3};
\draw[red] (153,225) node [anchor=north west][inner sep=0.75pt]   [align=left]  {4};
\draw[red] (195,90) node [anchor=north west][inner sep=0.75pt]   [align=left]  {5};
\draw[blue] (89,180) node [anchor=north west][inner sep=0.75pt]   [align=left]  {1};
\draw[blue] (137,145) node [anchor=north west][inner sep=0.75pt]   [align=left]  {2};
\draw[blue] (118,85) node [anchor=north west][inner sep=0.75pt]   [align=left]  {3};
\draw[blue] (58,85) node [anchor=north west][inner sep=0.75pt]   [align=left]  {4};
\draw[blue] (38,145) node [anchor=north west][inner sep=0.75pt]   [align=left]  {5};
\end{tikzpicture}
\end{displaymath}
\end{minipage}
\end{tabular}

\medskip

$F_5$ është  normalizuesi i   $C_5$ në $S_5$ dhe përkëmben pentagonin me të kundërtin e tij.  $D_5$ është nëngrupi  i $F_5$ që fikson pentagonin (grupi i simetrisë së pentagonit) dhe  $C_5$ nëngrupi i rrotullimeve.    Për shembull,  $F_5$   gjenerohet nga 
\begin{equation}
F_5=\< \sigma, \tau   \mid   \sigma^5=\tau^4=(\sigma\tau)^4=\sigma\sigma\tau\sigma^{-1} \tau^{-1} \>, 
\end{equation}
ku   $\sigma=(12345)$ and $\tau = (2453)$.   Pra në se  $G\le F_5$ atëherue   $G$ fikson

\begin{equation}\label{d_1}
\begin{split}
\d_1  = &     (\a_1-\a_2)^2 (\a_2-\a_3)^2 (\a_3-\a_4)^2 (\a_4-\a_5)^2 (\a_5-\a_1)^2 \\
&    -  (\a_1-\a_3)^2 (\a_3-\a_5)^2 (\a_5-\a_2)^2 (\a_2-\a_4)^2 (\a_4-\a_1)^2\\
\end{split}
\end{equation}
ku termi i parë (respektivisht termi i dytë) korrespondon me brinjët e pentagonit (respektivisht pentagonit të kundërt). 

Janë gjashtë nëngrupe  5-Sylow të $S_5$ që jepen si më poshtë
\[
\begin{split}
H_1 & = \<  (1,2,3,4,5)  \>  = \{   ( ),  (1,2,3,4,5), (1,3,5,2,4), (1,4,2,5,3), (1,5,4,3,2) \} \\
H_2 & =  \<  (1,2,3,5,4)  \>  =  \{   ( ),   (1,2,3,5,4), (1,3,4,2,5), (1,5,2,4,3), (1,4,5,3,2)  \} \\
H_3 & =  \< (1,2,4,5,3)   \>  =   \{   ( ),  (1,2,4,5,3), (1,4,3,2,5), (1,5,2,3,4), (1,3,5,4,2)  \} \\
H_4 & =   \< (1,2,4,3,5)   \>  =  \{   ( ),   (1,2,4,3,5), (1,4,5,2,3), (1,3,2,5,4),  (1,5,3,4,2)  \} \\
H_5 & =   \< (1,2,5,3,4)   \>  =   \{   ( ),  (1,2,5,3,4), (1,5,4,2,3), (1,3,2,4,5), (1,4,3,5,2)  \} \\
H_6 & =   \<  (1,3,4,5,2)  \>  =   \{   ( ),  (1,3,4,5,2),  (1,4,2,3,5), (1,5,2,4,3), (1,5,3,2,4)   \} \\
\end{split}
\]
Për të gjetur gjithë invariantet për çdo group $H_i$, $i=1, \dots , 6$ ne  konsiderojmë  polinomin  $f(x)$ si një formë binare $f(x, y)$.

\subsection{Format binare dhe invariantët e tyre}\label{invariants}
Një \textbf{formë binare} është një polinom homogjen $f(x, y)$ me gradë $d\geq 2$; shih \cite{b-sh} për një përshkrim të formave kuadratike në literaturën Shqip dhe \cite{Geyer, H} për format binare të gradës $d>2$.  Le të jetë $\F$ një fushë dhe $V_d (\F)$ hapësira e formave binare $f(x, y)$ me gradë   $d \geq 2$ dhe koefiçentë në $\F$.  $V_d (\F)$ është një hapësirë vektoriale mbi $\F$. Grupi    $\gl_2 ( \F)$ vepron mbi  $V_d (\F)$ si në vazhdim: 

\begin{align*}
\gl_2 ( k ) \times V_d  & \to V_d  \\
\left(    \begin{bmatrix} a & b \\ c & d \end{bmatrix}, f (x, y)  \right)    & \to   f (ax +by , cx + dy ) \\
\end{align*}
Ne do të shënojmë  $f (ax +by , cx + dy )$ me $f^M$.  Dy forma binare quhen \textbf{ekuivalente} kur janë në të njëjtën orbitë të këtij veprimi.

Le të marrim tani një formë binare   $f\in \F[x,y]$ me gradë $n$. Mbi mbylljen algjebrike $\bar \F$ të $\F$-së, kjo formë faktorizohet si më poshtë
\[
f(x,y) =  b_nx^n\ +\ b_{n-1}x^{n-1}y\ +  \dots  \ + \ b_0y^n =  (\b_1x - \g_1y)  \cdots  (\b_nx -\g_ny)    
\]
ku $\b_i,\g_i\in\bar \F$. Vërtetimi i këtij fakti mund të shihet lehtë si më poshtë.  Supozojmë se  $b_n\ne0$. (Në qoftë se  $b_n=0$ përdorim induksionin mbi $f/y$). Atëherë kemi 
\[
f(x, y) =  b_ny^n \left( \left(\frac{x}{y}\right)^n  + \frac{b_{n-1}}{b_n}   \left(\frac{x}{y}\right)^{n-1}   + \dots + \frac{b_{0}}{b_n}\right) = 
 b_n   y^n\ \left(      \frac{x}{y}-\a_1\right)   \cdots \left(\frac{x}{y}-\a_n   \right)
\]
për disa  $\a_i\in \bar \F$.   Vini re se çdo faktor mund të shkruhet si 
\[
(\b_i x - \g_i y)  = \det \begin{bmatrix} x & y \\ \g_i & \b_i \end{bmatrix} 
\]
Pra, 
\[
f =  \det
\begin{bmatrix}
x & y\\ 
\g_1  & \b_1
\end{bmatrix}
  \cdots\ 
\det
\begin{bmatrix} x & y \\
\g_n  & \b_n
\end{bmatrix}
\]
Veprimi i grupit  $\gl_2 (\F)$ mbi   $V_n (\F)$ për çdo  $M \in \gl_2(\F)$, 
\[
M  = 
\begin{bmatrix} a  & b \cr c  & d
\end{bmatrix}
\]
dërgon  formën binare $f$ tek
\[
f^M  :=  f (ax + by, cx + dy ) =   (\b_1'x - \g_1'y) \ \cdots\ (\b_n'x - \g_n'y).
\]
Pra  kemi
\[
(\g_i', \b_i') =  (\g_i, \b_i)\ M^{-1} \det(M)
\]
Versioni projektiv i  $\d_1$ është  $\tilde\d_1$, që formohet duke zëvendësuar  $\a_i-\a_j$ me 
\begin{equation}
D_{ij}  = \det 
\begin{bmatrix}
\g_i  & \b_i \\
\g_j  & \b_j \\
\end{bmatrix}
\end{equation}
në formulat që përkufizojnë  $\d_i$, $i=1, \ldots , 6$.
\begin{equation}
\begin{split}
\tilde\d_1 & =  D_{12}^2 D_{23}^2 D_{34}^2 D_{45}^2 D_{51}^2  -  D_{13}^2 D_{35}^2 D_{52}^2 D_{24}^2 D_{41}^2 \\
\end{split}
\end{equation}

Meqënëse  $S_5$ ka gjashtë  5-Sylow nëngrupe, shënojmë me  $\d_1, \dots,\d_6$ elementët përkatës të këtyre nëngrupeve.  
Këto elementë janë
\begin{equation}
\begin{split}
\tilde\d_2 & = D_{12}^2 D_{23}^2 D_{35}^2 D_{54}^2 D_{41}^2  -  D_{13 }^2 D_{34 }^2 D_{42 }^2 D_{25}^2 D_{51}^2  \\
\tilde\d_3 & =  D_{12}^2 D_{24}^2 D_{45}^2 D_{53}^2 D_{31}^2  -  D_{14}^2 D_{43}^2 D_{32}^2 D_{25}^2 D_{51}^2  \\ 
\tilde\d_4 & = D_{12}^2 D_{24}^2 D_{43}^2 D_{35}^2 D_{51}^2  -  D_{14}^2 D_{45}^2 D_{52}^2 D_{23}^2 D_{31}^2  \\  
\tilde\d_5 & = D_{12}^2 D_{25}^2 D_{53}^2 D_{34}^2 D_{41}^2  -  D_{15}^2 D_{54}^2 D_{42}^2 D_{23}^2 D_{31}^2  \\  
\tilde\d_6 & =  D_{13}^2 D_{34}^2 D_{45}^2 D_{52}^2 D_{21}^2  -  D_{14}^2 D_{42}^2 D_{23}^2 D_{35}^2 D_{51}^2  \\ 
\end{split}
\end{equation}

\begin{lem} 
 $\d_i^\sigma = \d_i$ dhe $\d_i^\tau = \d_i$ për $i=1, \dots , 6$.
\end{lem}

 $G$ përkëmben $\d_1, \dots,\d_6$. Nëse  $G$  është i konjuguar me një nëngrup të $F_5$, atëherë fikson një nga  $\d_1, \dots,\d_6$
 Pikërisht elementi    $\d_i$ që fiksohet duhet të jetë në   $\F$.  Pra, kusht i nevojshëm dhe i mjaftueshëm që një polinom i pathjeshtueshëm $f(x)$ të jetë i zgjidhshëm me radikale është që një nga $\d_i$ të jetë në $\F$, ose me fjalë të tjera polinomi 
\begin{equation}
g(x) =  (x-\d_1)  \cdots (x-\d_6)  \in \F[x]
\end{equation}
të ketë një rrënjë në $\F$, ku $\d_1$ jepet në \cref{d_1} dhe $\d_2, \ldots , \d_6$ si më poshtë:
\begin{equation}
\begin{split}
\d_2  = &     (\a_1-\a_2)^2 (\a_2-\a_3)^2 (\a_3-\a_5)^2 (\a_5-\a_4)^2 (\a_4-\a_1)^2 \\
&    -  (\a_1-\a_3)^2 (\a_3-\a_4)^2 (\a_4-\a_2)^2 (\a_2-\a_5)^2 (\a_5-\a_1)^2\\
\d_3  = &     (\a_1-\a_2)^2 (\a_2-\a_4)^2 (\a_4-\a_5)^2 (\a_5-\a_3)^2 (\a_3-\a_1)^2 \\
&    -  (\a_1-\a_4)^2 (\a_4-\a_3)^2 (\a_3-\a_2)^2 (\a_2-\a_5)^2 (\a_5-\a_1)^2\\
\d_4  = &     (\a_1-\a_2)^2 (\a_2-\a_4)^2 (\a_4-\a_3)^2 (\a_3-\a_5)^2 (\a_5-\a_1)^2 \\
&    -  (\a_1-\a_4)^2 (\a_4-\a_5)^2 (\a_5-\a_2)^2 (\a_2-\a_3)^2 (\a_3-\a_1)^2\\
\d_5  = &     (\a_1-\a_2)^2 (\a_2-\a_5)^2 (\a_5-\a_3)^2 (\a_3-\a_4)^2 (\a_4-\a_1)^2 \\
&    -  (\a_1-\a_5)^2 (\a_5-\a_4)^2 (\a_4-\a_2)^2 (\a_2-\a_3)^2 (\a_3-\a_1)^2\\
\d_6  = &     (\a_1-\a_3)^2 (\a_3-\a_4)^2 (\a_4-\a_5)^2 (\a_5-\a_2)^2 (\a_2-\a_1)^2 \\
&    -  (\a_1-\a_4)^2 (\a_4-\a_2)^2 (\a_2-\a_3)^2 (\a_3-\a_5)^2 (\a_5-\a_1)^2\\
\end{split}
\end{equation}

\begin{thm}
Në qoftë se  $G$ fikson një  $\d_i$ atëherë  $G$ është i konjuguar me një nëngrup të  $F_5$, kur  $\d_1, \dots,\d_6$ janë të ndryshme.
\end{thm}

\proof
Për të vërtetuar këtë mjafton të vërtetojmë se    $\d_1, \dots,\d_6$ janë dy nga dy të ndryshme kur  $\D_f\ne0$. Pra duhet të vërtetojmë se 
$ \D_f \neq 0 \implies \D_g\neq 0$. 
Duke përdorur algjebrën kompjutacionale ne mund të gjejmë diskriminantin $\D_g$ e $g(x)$ dhe shohim se $\frac {\D_g } {\D_f}$ është një polinom 
 si më poshtë
 \[
\D_g   =  \left( (\a_1-\a_2) ( \a_3-\a_4 ) ( \a_4-\a_5 )  ( \a_3-\a_ 5 \right) )^4 \cdot \D_f \cdot  I_2^2  \cdot I_3 \cdot I_4^2 \cdot I_6^2\\
 \]
ku $I_2, I_3, I_4$, dhe $I_6$ janë dhënë në \cref{eqs}.
Kjo përfundon vërtetimin. 

\qed

Koefiçentët e $g(x)$ janë funksione simetrike në $\a_1, \dots,\a_5$, prandaj janë shprehje polinomiale $c_0, \dots, c_4$. Qëllimi ynë është të gjejmë këto shprehje. Kjo do na japë një kriter për të kontrolluar nëse $f(x)=0$ është i zgjidhshëm me radikale.

\begin{lem} Le të jetë $s_r( x_1, \dots, x_6)$, $r=1, \dots,6$, polinomi simetrik elementar
\begin{equation}
s_r = \sum_{i_1<i_2< \dots <i_r} x_{i_1}x_{i_2}  \dots  x_{i_r}.
\end{equation}
Atëherë
\begin{equation}
d_r: = s_r ( \tilde\d_1, \dots, \tilde\d_6)
\end{equation}
 është një shprehje polinomiale homogjene në $b_0, \dots,b_5$ e gradës $4 r$. Këto polinomeve janë invariante nën veprimin e $\Sl_2(\F)$ në kuintikët binare: Për çdo $M\in \Sl_2(\F)$ kuintiku $f^M$ ka të shoqëruar të njëjtin $d_r$.
\end{lem}

\begin{figure}[hb]
\[
\xymatrix@R=4em@C=2em{
        &       &        &   \ar@{-}[dl]   \ar@{-}[dll]   \ar@{-}[dlll] \F(\a_1,\dots, \a_5)   \ar@{-}[dr]    \ar@{-}[drr]   \ar@{-}[drrr]    & &  &    \\
\F(\d_1)  \ar@{-}[drrr]  & \F(\d_2) \ar@{-}[drr] & \F(\d_3) \ar@{-}[dr] & &  \F(\d_4) \ar@{-}[dl] & \F(\d_5) \ar@{-}[dll] & \F(\d_6) \ar@{-}[dlll]     \\
&&&  \F(\d_1,\cdots, \d_6)  \ar@{-}[d] &&&\\
&&& \F (s_1, \dots , s_6) \ar@{-}[d] &&&\\
&&&   \F    &&&}
\]
\end{figure} 

\proof Për çdo  $\a_j:=\g_j/\b_j$ kemi
 \[
 \tilde\d_i=  (\b_1\cdots\b_5)^4  \d_i =  b_5^4 \cdot \d_i.
 \]
Pra 
 $d_r = b_5^{4r}    s_r(\d_1,, \dots, \d_6)$. Meqënëse  $ s_r(\d_1,, \dots, \d_6)$ janë shprehje polinomiale në  $c_j=\frac {b_j} {b_5}$, për   $j=0, \dots,4$, atëherë 
 $d_r$ \es funksion racional në  $b_0, \dots,b_5$, ku numëruesi \es një fuqi e  $b_5$. 
 Duke këmbyer rolet e  $x$ me $y$ kemi që edhe emëruesi \es një fuqi e  $b_0$. 
 Kjo do të thotë që është konstante, i.e., $d_r$ është një polinom në $b_0, \dots,b_5$. Nëse zevëndësojmë secilin $\b_j$ me $c\b_j$ për një skalar $\lambda$ atëherë  secili  $\tilde\d_i$ do të shumëzohet me $\lambda^4$, kështu $d_r$ do të shumëzohet me $\lambda^{4r}$. Prandaj $d_r$ është homogjen i gradës $4r$. Pjesa tjetër është e qartë.
 
\qed


%
\[
\begin{split}
I_2  = & \left( \a_1\a_3+\a_1\a_ 4-2 \a_1\a_5-2 \a_3\a_4+\a_3 \a_5+\a_4\a_5 \right)  \left( 2 \a_2\a_3-\a_ 2\a_4-\a_2\a_5-\a_3\a_4-\a_3\a_5+2 \a_4\a_5 \right)  \\
&  \left( \a_2\a_3-2 \a_2\a_4+\a_2\a_5+\a_3\a_4-2 \a_3\a_5+\a_4\a_5 \right)  \left( 2 \a_1\a_2-\a_1\a_4-\a_1\a_5-\a_2\a_4-\a_2\a_5+2 \a_4\a_5 \right) \\
&  \left( \a_1\a_3-2 \a_1\a_4+\a_1\a_5+\a_3\a_4-2 \a_3\a_5+\a_4\a_5 \right)  \left( \a_2\a_3+\a_2\a_4-2 \a_2\a_5-2 \a_3\a_4+\a_3\a_5+\a_4\a_5 \right)   \\
 & \left( 2 \a_1\a_3-\a_1\a_ 4-\a_1\a_5-\a_3\a_4-\a_3\a_5+2 \a_4\a_5 \right)  \left( 2 \a_1\a_2-\a_1\a_3-\a_1\a_5-\a_2\a_3-\a_2\a_5+2 \a_3\a_5 \right)   \\
 & \left( 2 \a_1\a_2-\a_1\a_3-\a_1\a_4-\a_2\a_3-\a_2\a_4+2 \a_3\a_4 \right)  
\end{split}
\]

\begin{small}
\[
\begin{split} 
 I_3 = 
& \left( \a_1\a_2\a_3-\a_1\a_2 \a_5-2 \a_1\a_3\a_4+\a_1 \a_3\a_5+\a_1\a_4^2-\a_2 \a_3\a_5-\a_2\a_4^2+2 \a_ 2\a_4\a_5+\a_3\a_4^2-\a_4^2\a_5 \right)^2 \\
&  \left( \a_1\a_2\a_4-\a_1\a_2\a_5-\a_1\a_3^2+2 \a_1\a_3\a_5-\a_1\a_4\a_5+\a_2\a_3^2-2 \a_2\a_3\a_4+\a_2\a_4\a_5+\a_3^2\a_4-\a_3^2\a_5 \right)^2 \\
& \left( \a_1\a_2\a_ 3-\a_1\a_2\a_4+\a_1\a_3\a_4-2 \a_1\a_3\a_5+\a_1\a_5^2-\a_2\a_3\a_4+2 \a_2\a_4\a_5-\a_2\a_5^2+\a_3\a_5^2-\a_4\a_5^2 \right)^2 \\
& \left( \a_1\a_2\a_4-\a_1\a_2\a_ 5+\a_1\a_3^2-2 \a_1\a_3\a_ 4+\a_1\a_4\a_5-\a_2\a_3^2+2 \a_2\a_3\a_5-\a_2\a_4\a_5+\a_3^2\a_4-\a_3^2\a_ 5 \right)^2 \\
& \left( \a_1\a_2\a_3-\a_1\a_2\a_ 4-\a_1\a_3\a_4+2 \a_1\a_4\a_5-\a_1\a_5^2+\a_2\a_3\a_4-2 \a_2\a_3\a_5+\a_2\a_5^2+\a_3\a_5^2-\a_4\a_5^2 \right)^2  \\
& \left( \a_1\a_2\a_3-\a_1\a_2\a_5-\a_1\a_3\a_5-\a_1\a_4^2+2 \a_1\a_4\a_5-2 \a_2\a_3\a_4+\a_2\a_3\a_ 5+\a_2\a_4^2+\a_3\a_4^2-\a_4^2\a_5 \right)^2 \\
& \left( \a_1 \a_2\a_3-2 \a_1\a_2\a_4+\a_1\a_2\a_5-\a_1\a_3\a_ 5+\a_1\a_4^2-\a_2\a_3\a_ 5+\a_2\a_4^2-\a_3\a_4^2+2 \a_3\a_4\a_5-\a_4^2\a_5 \right)^4  \\
&\left( \a_1\a_2\a_3+\a_1\a_2\a_4-2 \a_1 \a_2\a_5-\a_1\a_3\a_4+\a_1\a_5^2-\a_2\a_3\a_4+\a_ 2\a_5^2+2 \a_3\a_4\a_5-\a_ 3\a_5^2-\a_4\a_5^2 \right)^4 \\
& \left( 2 \a_1\a_2\a_3-\a_1\a_ 2\a_4-\a_1\a_2\a_5-\a_1\a_3^2+\a_1\a_4\a_5-\a_2\a_3^2+\a_2\a_4\a_5+\a_3^2\a_4+\a_3^2\a_5-2 \a_3\a_ 4\a_5 \right)^4 
\end{split}
\]

\[
\begin{split}
I_4 = &     ( \a_1^2\a_2^2-2 \a_1^2\a_2\a_5+\a_1^2\a_4^2-2 \a_1^2\a_4\a_5+2 \a_1^2\a_5^2-2 \a_1\a_2^2\a_4-2 \a_1\a_2\a_4^2 +8 \a_1\a_2\a_4\a_5-2 \a_1\a_2\a_5^2-2 \a_1\a_4\a_5^2        \\
& +2 \a_2^2\a_4^2-2 \a_2^2\a_ 4\a_5+\a_2^2\a_5^2-2 \a_2\a_4^2\a_5+\a_4^2\a_5^2 ) ( \a_1^2\a_2^2-2 \a_1^2\a_2\a_4+\a_1^2\a_3^2-2 \a_1^2\a_3\a_4+2 \a_1^2\a_4^2-2 \a_1\a_2^2\a_3  \\
& -2 \a_1\a_2\a_3^2  +8 \a_1\a_2\a_3\a_4-2 \a_1\a_2\a_4^2-2 \a_1\a_3\a_4^2+2 \a_2^2\a_3^2-2 \a_2^2\a_3\a_4+\a_2^2\a_4^2-2 \a_2\a_3^2\a_4+\a_3^2\a_4^2 )  ( \a_1^2\a_2^2 \\
& -2 \a_1^2\a_2\a_3+2 \a_1^2\a_3^2-2 \a_1^2\a_3\a_5+\a_1^2\a_5^2-2 \a_1\a_2^2\a_5-2 \a_1\a_2\a_3^2  +8 \a_1\a_ 2\a_3\a_5-2 \a_1\a_2\a_5^2-2 \a_1\a_3^2\a_5+\a_2^2\a_3^2 \\
& -2 \a_2^2\a_3\a_5+2 \a_2^2\a_5^2-2 \a_2\a_3\a_5^2+\a_3^2\a_5^2 ) ( \a_1^2\a_2^2-2 \a_1^2\a_2\a_3+2 \a_1^2\a_3^2-2 \a_1^2\a_3\a_4+\a_1^2\a_4^2-2 \a_1\a_2^2\a_4-2 \a_1\a_2\a_3^2  \\
& +8 \a_1\a_2\a_3\a_4-2 \a_1\a_2\a_4^2-2 \a_1\a_3^2\a_4+\a_2^2\a_3^2-2 \a_2^2\a_3\a_4+2 \a_2^2\a_4^2-2 \a_2\a_3\a_4^2+\a_3^2\a_4^2 ) ( \a_1^2\a_2^2-2 \a_1^2\a_2\a_4+2 \a_1^2\a_4^2\\
& -2 \a_1^2\a_4\a_5+\a_1^2\a_5^2-2 \a_1\a_2^2\a_5-2 \a_1\a_2\a_4^2  +8 \a_1\a_2\a_4\a_5-2 \a_1\a_2\a_5^2-2 \a_1\a_4^2\a_5+\a_2^2\a_4^2-2 \a_2^2\a_4\a_5+2 \a_2^2\a_5^2\\
& -2 \a_2\a_4\a_5^2 +\a_4^2\a_5^2 )( \a_1^2\a_2^2-2  \a_1^2\a_2\a_5+\a_1^2{\a_ 3}^2-2 \a_1^2\a_3\a_5+2 \a_1^2\a_5^2-2 \a_1\a_2^2\a_ 3-2 \a_1\a_2\a_3^2   +8 \a_1\a_2\a_3\a_5\\
& -2 \a_1\a_2\a_5^2-2 \a_1\a_3\a_5^2 +2 \a_2^2\a_3^2-2 \a_2^2\a_3\a_5+\a_2^2\a_5^2-2 \a_2\a_5\a_3^2+\a_3^2\a_5^2 ) 
\end{split}
\]
\end{small}

Për $I_6$ shihni \cref{i6}. 


Janë katër invariante kryesore të kuintikëve, të shënuara $J_4, J_8, J_{12},  J_{18}$, të gradës 4,8,12 dhe 18, të tillë që secili $\Sl(2,\F)$- polinom invariant në $b_0, \dots,b_5$ është polinom $J_4,J_8, J_{12}, J_{18}$ (shiko e.g. I. Schur, Vorlesungen ueber Invariantentheorie, Springer 1968). Për të përcaktuar $J_4,  J_8,  J_{12}$, na duhen disa madhësi ndihmëse
\begin{equation*}
\begin{split}
A   &= \frac 1 {100} \left(20 b_4-8 b_1 b_3+3 b_2^2\right), \\
B  &= \frac1 {100} \left( 100 b_5-12 b_1 b_4+2 b_2 b_3\right), \\
C  &= \frac 1 {100} \left(20 b_1 b_5-8 b_2 b_4+3 b_3^2\right)
\end{split}
\end{equation*}
dhe $D, E, F, G$ të përkufizuara nga
\begin{equation}
 \frac 1{1000} \det
 \begin{bmatrix}
 10u+2b_1v  &  2b_1u+b_2v  &  b_2u+b_3v  \\
 2b_1u+b_2v  &  b_2u+b_3v   &  b_3u+2b_4v  \\
 b_2u+b_3v  &  b_3u+2b_4v  &  2b_4u+10b_5v  \\
\end{bmatrix}
=    Du^3 + Eu^2v + F  uv^2 + Gv^3
\end{equation}
Atëherë, 
\begin{equation}
\begin{split}
J_4 & =   5^3 (B^2 -  4 AC)\\
 J_8 &  =   2^5 \cdot 5^6 \left[ 2A(3EG-F^2)-B(9DG-EF)+2C(3FD-E^2)\right]\\
 J_{12} &  =  -2^{10} \cdot 5^9 \cdot 3^{-1} \left[ 4(3EG-F^2)(3FD-E^2) -  (9DG-EF)^2 \right]
\end{split}
\end{equation}
Nga teoria klasike e invariantëve mund të shprehen $d_1, \ldots , d_6$ si shprehje në $J_4, J_8, J_{12}$, një rezultat i 
 shihni Berwick 1915;  shihni \cite{king-beyond-quartic}.
 
\begin{equation}
\begin{split}
d_1 & =    -10 J_4\\
d_2 & =    35 J_4^2 + 10  J_8\\
d_3 & =    -60 J_4^3 - 30 J_4 J_8 - 10  J_{12}\\
d_4 & =    55 J_4^4 + 30 J_4^2 J_8 + 25  J_8^2 + 50 J_4 J_{12}\\
d_5 & =    -26 J_4^5 - 10 J_4^3 J_8 - 44 J_4 J_8^2 - 59 J_4^2 J_{12}  - 14  J_8 J_{12}\\
d_6 & =    5J_4^6 + 20 J_4^2 J_8^2 + 20 J_4^3 J_{12} + 20 J_4 J_8 J_{12} + 25  J_{12}^2\\
\end{split}
\end{equation}


Pra koefiçentët e $g(x)$ tani jepen si polinome në koefiçentët e $f(x)$.  Rezultatet e mësipërme i përmbledhim më poshtë:

\begin{thm}\label{gal_5}
Le të jetë dhënë  $f(x)$ një polinom kuintik i pathjeshtueshëm  me koefiçentë në  $\F$ dhe  $d_1, \dots, d_6$ të përkufizuara si më lart në koefiçentët e $f(x)$. 
Atëherë $f(x)$ është i zgjidhshëm me radikale  atëherë dhe vetëm atëherë kur 
\begin{equation}
g(x)= x^6 + d_1 x^5 + \cdots d_5 x + d_6
\end{equation}
ka një rrënjë në $\F$.
\end{thm}

\bibliographystyle{alpha} 

\bibliography{paper}{}

\restoregeometry

\appendix

\begin{landscape}

\section{Invariants}

\begin{small}
\[
\begin{split}
I_6 = &  ( \a_1^4\a_2^2-2 \a_1^4\a_2\a_3+\a_1^4\a_3^2+\a_1^4\a_4^2-2 \a_1^4\a_4\a_ 5+\a_1^4\a_5^2-2 \a_1^3\a_2^2\a_4   -2 \a_1^3\a_2^2\a_5+4 \a_1^3\a_2\a_3\a_4+4 \a_1^3\a_2\a_3\a_5-2 \a_1^3\a_2\a_4^2+4 \a_1^3\a_2\a_4\a_5  -2 \a_1^3\a_2\a_5^2 \\
& -2 \a_1^3\a_3^2\a_4-2 \a_1^3\a_3^2\a_5-2 \a_1^3\a_ 3\a_4^2+4 \a_1^3\a_3\a_4\a_5-2 \a_1^3\a_3\a_5^2+2 \a_1^2\a_2^2\a_4^2+2 \a_1^2\a_2^2\a_4\a_5+2 \a_1^2\a_2^2\a_5^2+2 \a_1^2\a_2\a_3\a_4^2-16 \a_1^2\a_2\a_3\a_4\a_5+2 \a_1^2\a_2\a_3\a_5^2\\
& +2 \a_1^2\a_3^2\a_4^2  +2 \a_1^2\a_3^2\a_4\a_5+2 \a_1^2\a_3^2\a_5^2-2 \a_1\a_2^2\a_3\a_4^2+4 \a_1\a_2^2\a_3\a_4\a_5-2 \a_1\a_2^2\a_3\a_5^2-2 \a_1\a_2^2\a_4^2\a_5-2 \a_1\a_2^2\a_4\a_5^2-2 \a_1\a_2\a_3^2\a_4^2+4 \a_1\a_2\a_3^2\a_4\a_5\\
& -2 \a_1\a_2\a_3^2\a_5^2  +4 \a_1\a_2\a_3\a_4^2\a_5  +4 \a_1\a_2\a_3\a_4\a_5^2-2 \a_1\a_3^2\a_4^2\a_5-2 \a_1\a_3^2\a_4\a_5^2+\a_2^2\a_3^2{\a_ 4}^2-2 \a_2^2\a_3^2\a_4\a_ 5+\a_2^2\a_3^2\a_5^2+\a_2^2\a_4^2\a_5^2-2 \a_2\a_ 3\a_4^2\a_5^2+\a_3^2\a_4^2\a_5^2 ) \\   
& ( \a_1^2\a_2^4-2 \a_1^2\a_2^3\a_3-2 \a_1^2\a_2^3\a_5+2 \a_1^2\a_2^2\a_3^2+2 \a_1^2\a_2^2\a_3\a_5+2 \a_1^2\a_2^2\a_5^2-2 \a_1^2\a_2\a_3^2\a_4-2 \a_1^2\a_2\a_3^2\a_5+4 \a_1^2\a_2\a_3\a_4\a_5-2 \a_1^2\a_2\a_3\a_5^2-2 \a_1^2\a_2\a_4\a_5^2+\a_1^2\a_3^2\a_4^2 \\
& +\a_1^2\a_3^2\a_5^2-2 \a_1^2\a_ 3\a_4^2\a_5+\a_1^2\a_4^2\a_5^2-2 \a_1\a_2^4\a_4-2 \a_1\a_2^3\a_3^2+4 \a_1\a_2^3\a_3\a_4+4 \a_1\a_2^3\a_3\a_5+4 \a_1\a_2^3\a_4\a_5-2 \a_1\a_2^3\a_5^2+2 \a_1\a_2^2\a_3^2\a_4-16 \a_1\a_2^2\a_3\a_4\a_5 \\
& +2 \a_1\a_2^2\a_4\a_5^2-2 \a_1\a_2\a_3^2\a_4^2+4 \a_1\a_2\a_3^2\a_4\a_5+4 \a_1\a_2\a_3\a_4^2\a_5+4 \a_1\a_2\a_3\a_4\a_5^2-2 \a_1\a_2\a_4^2\a_5^2-2 \a_1\a_3^2\a_4\a_5^2+\a_2^4\a_3^2-2 \a_2^4\a_3\a_5+\a_2^4\a_4^2+\a_2^4\a_5^2\\
& -2 \a_2^3\a_3^2\a_4 -2 \a_2^3\a_3\a_4^2+4 \a_2^3\a_3\a_4\a_5-2 \a_2^3\a_4^2\a_5-2 \a_2^3\a_4\a_5^2+2 \a_2^2\a_3^2\a_4^2+2 \a_2^2\a_3\a_4^2\a_5+2 \a_2^2\a_4^2\a_5^2-2 \a_2\a_3^2\a_4^2\a_5-2 \a_2\a_3\a_4^2\a_5^2+\a_3^2\a_4^2\a_5^2 )  \\
&  \a_1^2\a_2^2\a_3^2-2 \a_1^2\a_2^2\a_3\a_4+\a_1^2\a_2^2\a_4^2-2 \a_1^2\a_2\a_3^2\a_5+4 \a_1^2\a_2\a_3\a_4\a_5-2 \a_1^2\a_2\a_4^2\a_5+\a_1^2\a_3^2\a_4^2-2 \a_1^2\a_3^2\a_4\a_5+2 \a_1^2\a_3^2\a_5^2-2 \a_1^2\a_3\a_4^2\a_5 \\
& +2 \a_1^2\a_3\a_4\a_5^2-2 \a_1^2\a_3\a_5^3+2 \a_1^2\a_4^2\a_5^2-2 \a_1^2\a_4\a_5^3+\a_1^2\a_5^4-2 \a_1\a_2^2\a_3^2\a_5+4 \a_1\a_2^2\a_3\a_4\a_5-2 \a_1\a_2^2\a_4^2\a_5-2 \a_1\a_2\a_3^2\a_4^2+4 \a_1\a_2\a_3^2\a_4\a_5  +2\a_1\a_2\a_3^2\a_5^2  \\
& +4 \a_1\a_2\a_3\a_4^2\a_5-16 \a_1\a_2\a_3\a_4\a_5^2+4 \a_1\a_2\a_3\a_5^3+2 \a_1\a_ 2\a_4^2\a_5^2+4 \a_1\a_2\a_4\a_5^3-2 \a_1\a_2\a_5^4-2 \a_1\a_3^2\a_5^3+4 \a_1\a_3\a_4\a_5^3-2 \a_1\a_4^2\a_5^3 +\a_2^2\a_3^2\a_4^2-2\a_2^2\a_3^2\a_4\a_5\\
&  +2 \a_2^2\a_3^2\a_5^2-2 \a_2^2\a_3\a_4^2\a_5+2 \a_2^2\a_3\a_4\a_5^2-2 \a_2^2\a_3\a_5^3+2 \a_2^2\a_4^2\a_5^2-2 \a_2^2\a_4\a_5^3+\a_2^2\a_5^4-2 \a_2\a_3^2\a_5^3+4 \a_2\a_3\a_4\a_5^3-2 \a_2\a_4^2\a_5^3+\a_3^2\a_5^4-2\a_3\a_4\a_5^4+\a_4^2\a_5^4 ) \\
&   ( \a_1^4\a_2^2-2 \a_1^4\a_2\a_4+\a_1^4\a_3^2-2 \a_1^4\a_3\a_5+\a_1^4{\a_ 4}^2+\a_1^4\a_5^2-2 \a_1^3\a_2^2\a_3-2 \a_1^3\a_2^2\a_5-2 \a_1^3\a_2\a_3^2+4 \a_1^3\a_2\a_3\a_4+4 \a_1^3\a_2\a_3\a_5+4 \a_1^3\a_2\a_4\a_5-2 \a_1^3\a_2\a_5^2 \\
& -2 \a_1^3\a_3^2\a_4-2 \a_1^3\a_3\a_4^2+4 \a_1^3\a_3\a_4\a_5-2 \a_1^3{\a_ 4}^2\a_5-2 \a_1^3\a_4\a_5^2+2 \a_1^2\a_2^2\a_3^2+2 \a_1^2\a_2^2\a_3\a_5+2 \a_1^2\a_2^2\a_5^2+2 \a_1^2\a_ 2\a_3^2\a_4-16 \a_1^2\a_2\a_3\a_4\a_5+2 \a_1^2\a_2\a_4\a_5^2\\
& +2 \a_1^2\a_3^2\a_4^2+2 \a_1^2\a_3\a_4^2\a_5+2 \a_1^2\a_4^2\a_5^2-2 \a_1\a_2^2\a_3^2\a_4-2 \a_1\a_2^2\a_3^2\a_5+4 \a_1\a_2^2\a_3\a_4\a_5-2 \a_1\a_2^2\a_3\a_5^2-2 \a_1\a_2^2\a_4\a_5^2-2 \a_1\a_2\a_3^2\a_4^2+4 \a_1\a_2\a_3^2\a_4\a_5\\
& +4 \a_1\a_2\a_3\a_4^2\a_5+4 \a_1\a_2\a_3\a_4\a_5^2-2 \a_1\a_2\a_4^2\a_5^2-2 \a_1\a_3^2\a_4^2\a_5-2 \a_1\a_3{\a_ 4}^2\a_5^2+\a_2^2\a_3^2{\a_ 4}^2+\a_2^2\a_3^2\a_5^2-2 \a_2^2\a_3\a_4^2\a_5+\a_2^2\a_4^2\a_5^2-2 \a_2{\a_ 3}^2\a_4\a_5^2+\a_3^2\a_4^2\a_5^2  )  \\
&  ( \a_1^2\a_2^4-2 \a_1^2\a_2^3\a_3-2 \a_1^2\a_2^3\a_4+2 \a_1^2\a_2^2\a_3^2+2 \a_1^2\a_2^2\a_3\a_4+2 \a_1^2\a_2^2\a_4^2-2 \a_1^2\a_2\a_3^2\a_4-2 \a_1^2\a_2\a_3^2\a_5-2 \a_1^2\a_2\a_3\a_4^2+4 \a_1^2\a_2\a_3\a_4\a_5-2 \a_1^2\a_2\a_4^2\a_5\\
& +\a_1^2\a_3^2\a_4^2+\a_1^2\a_3^2\a_5^2-2 \a_1^2\a_ 3\a_4\a_5^2+\a_1^2\a_4^2\a_5^2-2 \a_1\a_2^4\a_5-2 \a_1\a_2^3\a_3^2+4 \a_1\a_2^3\a_3\a_4+4 \a_1\a_2^3\a_3\a_5-2 \a_1\a_2^3\a_4^2+4 \a_1\a_2^3\a_4\a_5+2 \a_1\a_2^2\a_3^2\a_5-16 \a_1\a_2^2\a_3\a_4\a_5\\
& +2 \a_1\a_2^2\a_4^2\a_5+4 \a_1\a_2\a_3^2\a_4\a_5-2 \a_1\a_2\a_3^2\a_5^2+4 \a_1\a_2\a_3\a_4^2\a_5+4 \a_1\a_2\a_3\a_4\a_5^2-2 \a_1\a_2\a_4^2\a_5^2-2 \a_1\a_3^2\a_4^2\a_5+\a_2^4\a_3^2-2 \a_2^4\a_3\a_4+\a_2^4\a_4^2+\a_2^4\a_5^2-2 \a_2^3\a_3^2\a_5\\
&+4 \a_2^3\a_3\a_4\a_5-2 \a_2^3\a_3\a_5^2-2 \a_2^3\a_4^2\a_5-2 \a_2^3\a_4\a_5^2+2 \a_2^2\a_3^2\a_5^2+2 \a_2^2\a_3\a_4\a_5^2+2 \a_2^2\a_4^2\a_5^2-2 \a_2\a_3^2\a_4\a_5^2-2 \a_2\a_3{\a_ 4}^2\a_5^2+\a_3^2\a_4^2\a_5^2  ) \\
\end{split}
\]

\[
\begin{split}
&  ( \a_1^2\a_2^2\a_4^2-2 \a_1^2\a_2^2\a_4\a_5+\a_1^2\a_2^2\a_5^2-2 \a_1^2\a_2\a_3\a_4^2+4 \a_1^2\a_2\a_3\a_4\a_5-2 \a_1^2\a_2\a_3\a_5^2+\a_1^2\a_3^4-2 \a_1^2\a_3^3\a_4-2 \a_1^2\a_3^3\a_5+2 \a_1^2\a_3^2\a_4^2+2 \a_1^2\a_3^2\a_4\a_5\\
& +2 \a_1^2\a_3^2\a_5^2-2 \a_1^2\a_3\a_4^2\a_5-2 \a_1^2\a_3\a_4\a_5^2+\a_1^2\a_4^2\a_5^2-2 \a_1\a_2^2\a_3\a_4^2+4 \a_1\a_2^2\a_3\a_4\a_ 5-2 \a_1\a_2^2\a_3\a_5^2-2 \a_1\a_2\a_3^4+4 \a_1\a_2\a_3^3\a_4+4 \a_1\a_2{\a_ 3}^3\a_5+2 \a_1\a_2\a_3^2\a_4^2\\
& -16 \a_1\a_2\a_3^2\a_4\a_5+2 \a_1\a_2\a_3^2\a_5^2+4 \a_1\a_2\a_3\a_4^2\a_5+4 \a_1\a_2\a_3\a_4\a_5^2-2 \a_1\a_2\a_4^2\a_5^2-2 \a_1\a_3^3\a_4^2+4 \a_1\a_3^3\a_4\a_5-2 \a_1\a_3^3\a_5^2+\a_2^2\a_3^4-2 \a_2^2\a_3^3\a_4-2 \a_2^2\a_3^3\a_5\\
& +2 \a_2^2{\a_ 3}^2\a_4^2+2 \a_2^2\a_3^2\a_4\a_5+2 \a_2^2\a_3^2\a_5^2-2 \a_2^2\a_3\a_4^2\a_5-2 \a_2^2\a_3\a_4\a_5^2+\a_2^2\a_4^2\a_5^2-2 \a_2\a_3^3\a_4^2+4 \a_2\a_3^3\a_4\a_5-2 \a_2\a_3^3\a_5^2+\a_3^4\a_4^2-2 \a_3^4\a_4\a_5+\a_3^4\a_5^2  ) \\
&  ( \a_1^4\a_2^2-2 \a_1^4\a_2\a_5+\a_1^4{\a_ 3}^2-2 \a_1^4\a_3\a_4+\a_1^4\a_4^2+\a_1^4\a_5^2-2 \a_1^3\a_2^2\a_3-2 \a_1^3\a_2^2\a_4-2 \a_1^3\a_2\a_3^2+4 \a_1^3\a_2\a_3\a_4+4 \a_1^3\a_2\a_3\a_5-2 \a_1^3\a_2\a_4^2+4 \a_1^3\a_2\a_4\a_5 \\
& -2 \a_1^3\a_3^2\a_5+4 \a_1^3\a_3\a_4\a_5-2 \a_1^3\a_3\a_5^2-2 \a_1^3\a_4^2\a_5-2 \a_1^3\a_4\a_5^2+2 \a_1^2\a_2^2\a_3^2+2 \a_1^2\a_2^2\a_3\a_4+2 \a_1^2\a_2^2\a_4^2+2 \a_1^2\a_2\a_3^2\a_5-16 \a_1^2\a_2\a_3\a_4\a_5+2 \a_1^2\a_2\a_4^2\a_5\\
& +2 \a_1^2\a_3^2\a_5^2+2 \a_1^2\a_3\a_4\a_5^2+2 \a_1^2\a_4^2\a_5^2-2 \a_1\a_2^2\a_3^2\a_4-2 \a_1\a_2^2\a_3^2\a_5-2 \a_1\a_2^2\a_3\a_4^2+4 \a_1\a_2^2\a_3\a_4\a_5-2 \a_1\a_2^2\a_4^2\a_5+4 \a_1\a_2\a_3^2\a_4\a_5-2 \a_1\a_2\a_3^2\a_5^2\\
& +4 \a_1\a_2\a_3\a_4^2\a_5+4 \a_1\a_2\a_3\a_4\a_5^2-2 \a_1\a_2\a_4^2\a_5^2-2 \a_1\a_3^2\a_4\a_5^2-2 \a_1\a_3\a_4^2\a_5^2+\a_2^2{\a_ 3}^2\a_4^2+\a_2^2\a_3^2\a_5^2-2 \a_2^2\a_3\a_4\a_5^2+\a_2^2\a_4^2\a_5^2-2 \a_2\a_3^2\a_4^2\a_5+{\a_ 3}^2\a_4^2\a_5^2  ) \\
&  ( \a_1^2\a_2^2\a_3^2-2 \a_1^2\a_2^2\a_3\a_5+\a_1^2\a_2^2\a_5^2-2 \a_1^2\a_2\a_3^2\a_4+4 \a_1^2\a_2\a_3\a_4\a_5-2 \a_1^2\a_2\a_4\a_5^2+2 \a_1^2\a_3^2\a_4^2-2 \a_1^2\a_3^2\a_4\a_5+\a_1^2\a_3^2\a_5^2-2 \a_1^2\a_3\a_4^3+2 \a_1^2\a_3\a_4^2\a_5\\
& -2 \a_1^2\a_3\a_4\a_5^2+\a_1^2\a_4^4-2 \a_1^2\a_4^3\a_5+2 \a_1^2\a_4^2\a_5^2-2 \a_1\a_2^2\a_3^2\a_4+4 \a_1\a_2^2\a_3\a_4\a_5-2 \a_1\a_2^2\a_4\a_5^2+2 \a_1\a_2{\a_ 3}^2\a_4^2+4 \a_1\a_2\a_3^2\a_4\a_5-2 \a_1\a_2\a_3^2\a_5^2+4 \a_1\a_2\a_3{\a_ 4}^3\\
& -16 \a_1\a_2\a_3\a_4^2\a_5+4 \a_1\a_2\a_3\a_4\a_5^2-2 \a_1\a_2\a_4^4+4 \a_1\a_2\a_4^3\a_5+2 \a_1\a_2\a_4^2\a_5^2-2 \a_1\a_3^2\a_4^3+4 \a_1\a_3\a_4^3\a_5-2 \a_1\a_4^3\a_5^2+2 \a_2^2\a_3^2\a_4^2-2 \a_2^2\a_3^2\a_4\a_5 \\
& +\a_2^2\a_3^2\a_5^2-2 \a_2^2\a_3\a_4^3+2 \a_2^2\a_3\a_4^2\a_5-2 \a_2^2\a_3\a_4\a_5^2+\a_2^2\a_4^4-2 \a_2^2\a_4^3\a_5+2 \a_2^2\a_4^2\a_5^2-2 \a_2\a_3^2\a_4^3+4 \a_2\a_3\a_4^3\a_5-2 \a_2\a_4^3\a_5^2+{\a_ 3}^2\a_4^4-2 \a_3\a_4^4\a_5+\a_4^4\a_5^2  ) \\
&  ( \a_1^2\a_2^4-2 \a_1^2\a_2^3\a_4-2 \a_1^2\a_2^3\a_5+2 \a_1^2\a_2^2\a_4^2+2 \a_1^2\a_2^2\a_4\a_5+2 \a_1^2\a_2^2\a_5^2-2 \a_1^2\a_2\a_3\a_4^2+4 \a_1^2\a_2\a_3\a_4\a_5-2 \a_1^2\a_2\a_3\a_5^2-2 \a_1^2\a_2\a_4^2\a_5-2 \a_1^2\a_2\a_4\a_5^2 \\
& +\a_1^2\a_3^2\a_4^2-2 \a_1^2\a_3^2\a_4\a_5+\a_1^2\a_3^2\a_5^2+\a_1^2\a_4^2\a_5^2-2 \a_1\a_2^4\a_3+4 \a_1\a_2^3\a_3\a_4+4 \a_1\a_2^3\a_3\a_5-2 \a_1\a_2^3\a_4^2+4 \a_1\a_2^3\a_4\a_5-2 \a_1\a_2^3\a_5^2+2 \a_1\a_2^2\a_3\a_4^2-16 \a_1\a_2^2\a_3\a_4\a_ 5 \\
& +2 \a_1\a_2^2\a_3\a_5^2-2 \a_1\a_2\a_3^2\a_4^2+4 \a_1\a_2\a_3^2\a_4\a_5-2 \a_1\a_2\a_3^2\a_5^2+4 \a_1\a_2\a_3\a_4^2\a_5+4 \a_1\a_2\a_3\a_4\a_5^2-2 \a_1\a_3\a_4^2\a_5^2+\a_2^4\a_3^2+\a_2^4\a_4^2-2 \a_2^4\a_4\a_5+\a_2^4\a_5^2-2 \a_2^3\a_3^2\a_4\\
& -2 \a_2^3\a_3^2\a_5-2 \a_2^3\a_3\a_4^2+4 \a_2^3\a_3\a_4\a_5-2 \a_2^3\a_3\a_5^2+2 \a_2^2\a_3^2\a_4^2+2 \a_2^2\a_3^2\a_4\a_5+2 \a_2^2\a_3^2\a_5^2-2 \a_2\a_3^2\a_4^2\a_5-2 \a_2\a_3^2\a_4\a_5^2+\a_3^2\a_4^2\a_5^2  ) 
\end{split}
\]
\end{small}
\end{landscape}

\end{document}